\documentclass[12pt,twoside]{article}

\usepackage[cp1251]{inputenc}    
\usepackage[english, russian]{babel}    
\usepackage{imm16_4}
\usepackage{amsmath}            
\usepackage{amsfonts,amssymb}   
\usepackage{graphicx}          
\usepackage{cite}               

\newcommand{\doc}{{\proof}}
\newtheorem{zam}{\mbox{З а м е ч а н и е}\ \ \\ }
\newtheorem{example}{\mbox{П р и м е р}\ \ \\ }

\newcommand{\rav}{\stackrel{\triangle}{=}}
\newcommand{\ph}{\varphi}
\newcommand{\la}{\langle}
\newcommand{\ra}{\rangle}
\newcommand{\epsi}{\varepsilon}
\newcommand{\beg}{\varkappa}

\newcommand{\rref}[1]{$(\ref{#1})$}
\newcommand{\mm}[1]{{\bf{#1}}}
\newcommand{\ct}[1]{{\mathcal{#1}}}
\newcommand{\fr}[1]{{\mathfrak{#1}}}
\newcommand{\td}[1]{\widetilde{#1}}
\newcommand{\bo}{{\hfill {$\Box$}}}

\newcommand{\leqref}[1]{\stackrel{(\ref{#1})}{\leq}}

\usepackage[active]{srcltx}

\newcommand{\Leftclt}{Д.\,В.\,Хлопин}     
\newcommand{\Rightclt}{Условия оптимальности при устойчивости сопряжённой переменной}     
\newcommand{\UDK}{517.977.5}              
\begin{document}

\Mainclt   



\begin{Titul}
{\bf{НЕОБХОДИМЫЕ УСЛОВИЯ ОПТИМАЛЬНОСТИ В СЛУЧАЕ УСТОЙЧИВОСТИ
СОПРЯЖЕННОЙ ПЕРЕМЕННОЙ \\ ДЛЯ ЗАДАЧ НА БЕСКОНЕЧНОМ
ПРОМЕЖУТКЕ}\footnote{Работа частично поддержана программой президиума
РАН "Математическая теория управления".}
\\[2ex]
Д.\,В.\,Хлопин}\\[3ex]
\end{Titul}

\begin{Anot}
 Данная работа посвящена исследованию необходимых условий
  оптимальности для задач управления  на бесконечном промежутке
  времени. В первую очередь исследуются условия, при которых
  для оптимальности траектории необходимо, чтобы
  при некотором решении соотношений принципа максимума
  сопряженная переменная стремилась к нулю на бесконечности.
  Показано, что для необходимости этого достаточно, чтобы
  сопряженная переменная, как компонента решения принципа максимума,
  была устойчива по Ляпунову.
  Предложена также целая серия более тонких условий на сопряженную переменную,
  показано, как свести проверку их необходимости к
  проверке устойчивости специально выписанной  системы уравнений.

 В качестве применения у
  полученной в \cite{kr_asD},\cite{kr_as}
  формуле Коши для сопряженной переменной вдоль оптимальной
  траектории
    найдены существенно более широкие условия применимости.
  Показано, что в случае  конечности и непрерывной зависимости от
    начальных условий, используемых в формуле несобственных
    интегралов, каждое оптимальное управление обязано иметь
    единственное решение принципа максимума, удовлетворяющее
    некоторому условию трансверсальности, и это --- решение,
    полученное при помощи  формулы Коши.
\end{Anot}

\sloppy



\section*{Введение} \setcounter{equation}{0}

   Принцип максимума Понтрягина для
   задач управления на бесконечном промежутке был сформулирован
   уже в в классической монографии \cite{ppp}, однако эти соотношения
      не полны, и выделяют, вообще говоря, слишком широкое
     семейство подозрительных на экстремум траекторий.
     Дополнительно требуется условие на бесконечности,
     предложено      достаточно много вариантов таких краевых условий,
     однако как показано, например, в \cite{Halkin}, \cite[Example 10.2]{norv},
     \cite[примеры  6.1-6.6]{kr_as}, эти условия могут оказаться
       несовместными с соотношениями принципом максимума, а могут тривиально
       из них следовать.

       Данная работа посвящена поиску условий на
        сопряженную переменную, позволяющих выделять для всякого
        оптимального управления не более одного решения
        соотношений принципа максимума. Само доказательство
        необходимости соответствующих условий сводится в работе к
        проверке устойчивости по части переменных решений системы
        принципа максимума. Показано, что устойчивость произведения сопряженной переменной
        на матричнозначную функцию гарантирует необходимость соответствующего
        условия типа трансверсальности. Оказывается, что если в качестве матричнозначной функции
        взять матричную экспоненту вдоль линеаризованного оптимального решения,
        то соответствующее условие выделяет единственную экстремаль.
        Проверка необходимости этого условия сводится к проверке
        непрерывной зависимости от параметра несобственного интеграла
        формулы Коши.

  \section{Определения и обозначения}
  \setcounter{equation}{0}

 В качестве промежутка времени будем рассматривать
  $\mm{T}\rav\{t\in\mm{R}\,|\,t\geq 0\}.$ 
 Через $E$ (возможно снабженное какими-либо
 индексами) будем обозначать какие-либо вспомогательные конечномерные евклидовы
 пространства, норму в них будем обозначать $||\cdot||_E$.
 Через $(comp)(E)$ будем обозначать пространство всевозможных
 компактных подмножеств множества $E$, снабдим это пространство
 метрикой Хаусдорфа.

  На множествах непрерывных на всем $\mm{T}$  функций будет
  рассматриваться как топология равномерной на всем $\mm{T}$
  сходимости, так и компактно-открытая топология ---
  топология сходимости равномерной на всяком компакте.
  (например $C(\mm{T},E)$ и $C_{loc}(\mm{T},E)$), первую
  из них снабдим нормой $||\cdot||_C$ топологии равномерной
  сходимости. Аналогично множества всех измеримых по Борелю
  функцмй будут снабжаться как той так и другой топологией
  (например $L^1([0,1],\mm{T})$ и $L^1_{loc}(\mm{T},\mm{T})$).

 Всюду далее, для всякой скалярной (или векторной) измеримой функции $a$
 интеграл $\int_{\mm{T}}a(t)dt$ будем понимать в несобственном смысле
 --- как предел $\int_{[0,T]}a(t)dt$ при $T\to\infty.$ Аналогично
 понимается интеграл по бесконечному промежутку, например $[T,\infty\ra.$

 Для всякого  подмножества $A$ топологического пространства под
  $cl\, A$  будем  понимать его замыкание.

 Обозначим через $\Omega$ семейство тех функций $\omega\in C(\mm{T},\mm{T}),$
 для которых $\lim_{t\to\infty}\omega(t)=0.$

   В качестве фазового пространства исходной управляемой системы будем использовать некоторое
   конечномерное
 метрическое пространство $\mm{X}\rav\mm{R}^m$. Единичный шар  в нём
 обозначим через $\mm{D}$.

 Пусть дано также конечномерное евклидовое пространство $\mm{U}$ и
   задано  отображение  $U:\mm{T}\mapsto (comp)(\mm{U})$,
  для которого выполнено:

  Условие ${\bf{(u)}}:$
  $U$\ ---\ ограниченное на всяком компактном множестве отображение,
  график $Gr\, U$ которого\ ---\  замкнутое борелевское множество в  $\mm{T}\times\mm{U}.$

 Будем говорить, что функция $a:\mm{T}\times E'\times \mm{U}\mapsto E''$
 удовлетворяет условиям Каратеодори, если
 1) функция $a(\cdot,y,u):\mm{T}\mapsto E''$ измерима для всех
 $(y,t,u)\in\mm{X}\times Gr {U},$
 2) функция $a(t,\cdot,\cdot):E'\times U(t)\mapsto E''$ непрерывна при
 всех $t\in\mm{T}.$

 Будем говорить, что функция $a:\mm{T}\times E'\times \mm{U}\mapsto E''$
 локально
 липшицева, если
  для любого компакта $K\in(comp)(E\times Gr\, U)$ найдется такая
  суммируемая на всяком компакте функция $L_K^a:\mm{T}\mapsto \mm{T},$ что для всех
  $(x',t,u),(x'',t,u)\in K$ выполнено
  $||a(t,x',u)-a(t,x'',u) ||_{E''}\leq L_K^a(t)||x'-x''||_{E'}.$

 Будем говорить, что функция $a:\mm{T}\times E'\times \mm{U}\mapsto E''$
 имеет измеримую мажоранту, если
  для любого компакта $K\in(comp)(E\times Gr\, U)$ найдётся такая
  суммируемая  на всяком компакте функция $M_K^a:\mm{T}\mapsto \mm{T},$ что для всех
  $(x,t,u)\in K$ выполнено
  $||a(t,x,u)||_{E''}\leq M_K^a(t).$

 Будем говорить, что для функции  $a:\mm{T}\times E'\times \mm{U}\mapsto E''$
 выполнено условие продолжимости
 решений на $\mm{T}$, если
 выполнено
 условие подлинейного по $x$ роста (для этого, в свою очередь достаточно, чтобы она
 имела измеримую мажоранту
 и была липшицева с независящей от $K$ функцией $L_K^a$; см.\cite[1.4.6]{tovst1}).


   Пусть дана управляемая система
   \begin{equation}
   \label{sys}
    \dot{x}=f(t,x,u),\ x(0)=0,\ t\in\mm{T},\ x\in\mm{X}, u\in U(t).
   \end{equation}

 Всюду далее предполагается, что выполнено:

 Условие ${\bf{(f)}}:$ отображение $f:\mm{T}\times\mm{X}\times \mm{U}\mapsto\mm{X}$
 --- отображение Каратеодори, локально липшицевое, имеющее измеримую
 мажоранту и удовлетворяющее условию продолжимости.

   Под множеством всевозможных допустимых управлений $\fr{U}$ будем понимать
  множество всех измеримых селекторов многозначного отображения $U$.
  Топологию на $\fr{U}$ зададим в силу вложения
  $\fr{U}\subset B_{loc}(\mm{T},\mm{U}).$
 Теперь всякому  $u\in\fr{U}$  можно сопоставить решение уравнения \rref{sys}, это решение единственно в силу локальной липшицевости $f$,
 в силу условия
 продолжимости оно может быть продолжено на все $\mm{T},$ обозначим его через
  $\ph[u].$ Теперь в силу \cite[теорема 1.1.6]{f} имеет место на каждом компакте
   непрерывная зависимость
  решения от программного управления, таким образом отображение
  $\ph:\fr{U}\to C_{loc}(\mm{T},\mm{X})$ непрерывно.
Предположим, что всегда далее выполнено:

 Условие ${\bf{(g)}}:$  отображение $g:\mm{T}\times\mm{X}\times \mm{U}\mapsto\mm{R}$
 --- локально липшицевое и имеющее измеримую
 мажоранту отображение Каратеодори, а кроме того
 для некоторой функции $\omega\in\Omega$
 при любых $u\in \fr{U}$ для всех $T\in\mm{T}$ выполнено
 $\displaystyle\int_{[T,\infty\ra}\Big|g\big(t,\ph[u](t),u(t)\big)\Big|dt \leq \omega(T).$

 Всюду далее исследуем задачу максимизации  на траекториях системы \rref{sys}
 функционала
 \begin{equation}
   \label{opt}
   J(u)\rav\int_{\mm{T}} g\big(t,\ph[u](t),u(t)\big)dt\to\max.
   \end{equation}
 В силу условия ${\bf{(g)}}$ функционал $J:\fr{U}\to\mm{R}$ непрерывен.


\section{Множество обобщенных управлений $\td{\fr{U}}$}
\setcounter{equation}{0}

  Обозначим для всех $u\in\mm{U}$ через $\td\delta_u$
   вероятностную меру, сосредоточенную в точке
   $u\in\mm{U}$. Если некоторое отображение $\phi:\fr{U}\to C_{loc}(\mm{T},E)$
  таково, что для всякого $t\in\mm{T}$ из $u'|_{[0,t]}=u''|_{[0,t]}$
  следует $\phi(u')|_{[0,t]}=\phi(u'')|_{[0,t]},$ то назовём
  такое отображение неупреждающим.

  Определим $\td{\fr{U}}_n$ как семейство всех слабо измеримых отображений
  $\mu$ из $[0,n]$
  в множество вероятностных мер Радона над $ \mm{U}$
   таких,
  что $\int_{U(t)} \eta(t)(du)=1$ для почти всех $t\in[0,n].$
    Оснастим это множество топологией *-слабой
  сходимости,
   полученное топологическое пространство --- компакт
   \cite[IV.3.11]{va},
   а множество
    $\fr{U}_n\rav\{u|_{[0,n]}\,|\,u\in\fr{U}\}$  всюду плотно вкладывается
   в $\td{\fr{U}}_n$ \cite[IV.3.10]{va} отображением $u\mapsto\td{\delta}\circ u$.

  Введем теперь $\td{\fr{U}}$\ ---
  семейство таких отображений  $\eta$ (из $\mm{T}$
  в множество вероятностных мер Радона над $ \mm{U}$),
  что  $\eta|_{[0,n]}\in \Pi([0,n], U)$ для всякого  $n\in\mm{N}$.
   Определим для всех $n\in\mm{N}$ сквозные проекции $\td\pi_n:\td{\fr{U}}\to \td{\fr{U}}_n$
   правилом $\td{\pi}_n(\eta)\rav\eta|_{[0,n]}$ для всех $\eta\in\td{\fr{U}}.$
  В качестве топологии на $\td{\fr{U}}$\ возьмем слабейшую топологию, в которой
  все сквозные проекции непрерывны; в частности
  множество $\ct{A}\subset\td{\fr{U}}$\ замкнуто, если и только если
  для всех $n\in\mm{N}$ образ $\td{\pi}_n(\ct{A})$ замкнут в
  $\td{\fr{U}}_n$.

   Каждой  $\eta\in \td{\fr{U}}$ можно сопоставить
  $\td\ph[\eta]\in C_{loc}(\mm{T},\mm{X})$ как
  решение задачи Коши
   $$\dot{x}=\int_{U(t)} f(\tau,x(\tau),u)\,\eta(t)(du),\ x(0)=0,$$
  и функционал $\td{J}(\eta)=\int_{\mm{T}}\int_{U(t)} g(\tau,\td\ph[\eta](\tau),u)\,\eta(t)(du)\,dt.$

   Множество $\td{\fr{U}}$ далее
   будем называть множеством обобщенных
   управлений. Введение этого множества оправдывает
   следующее предложение:

  \begin{prp}
  \label{spectr}
     Пусть выполнены ${\bf{(u),(f),(g)}}$, тогда

     1) пространство $\td{\fr{U}}$ --- компакт, а  $\td\delta({\fr{U}})$ --- всюду
     плотное его подмножество;

     2) всякому неупреждающему непрерывному отображению $\phi:\fr{U}\to C_{loc}(\mm{T},E)$ можно
     единственным образом сопоставить такое неупреждающее непрерывное отображение
     $\td\phi:\fr{U}\to C_{loc}(\mm{T},E),$ что $\phi=\td\phi\circ\td\delta,$
     при этом
     $cl\, \phi(\fr{U})=\td\phi(\td{\fr{U}});$

     3)  отображения $\td\ph,\td{J}$ непрерывны;

     4)    существует такой элемент $\eta^0\in\td{\fr{U}},$ что
\begin{equation}
   \label{supmax}
       \sup_{u\in\fr{U}}J(u)=\max_{\eta\in\td{\fr{U}}}\td{J}(\eta)=\td{J}(\eta^0).
   \end{equation}
\end{prp}
\doc
 Введем для краткости $\td\Pi\rav\prod_{n\in\mm{N}} \td{\fr{U}}_n$,
  снабдим произведение тихоновской топологией.
  Поскольку проекции $\td{\pi}_{n}$ непрерывны,
  то из \cite[II.1.3]{phillvv} будет непрерывно отображение
 $\td{\Delta}:\td{\fr{U}}\to\td\Pi,$ определенное по правилу:
 $\td{\Delta}(\eta)\rav\big(\td{\pi}_n(\eta)\big)_{n\in\mm{N}}$
  для всех
 $\eta\in\td{\fr{U}}.$
     Но отображение $\td{\Delta}$ инъективно, тогда
     существует обратное $\td{\Delta}^{-1}$.  Для проверки непрерывности $\td{\Delta}^{-1}$  достаточно
     проверить
     непрерывность суперпозиции $\td{\pi}_n\circ\td{\Delta}^{-1},$
      то есть операции проектирования из $\td{\Pi}$ на
     $n$-ю компоненту $\td{\Pi}$, что тривиально.

  Пусть  $n,k\in\mm{N}, (n>k).$  Тогда пространство $\td{\fr{U}}_n$
  вкладывается
  в $\td{\fr{U}}_k$  отображением
  $\td{\pi}^{n}_{k}(\eta)\rav\eta|_{[0,k]}$ для всех $\eta\in \td{\fr{U}}_n.$
  Поскольку при этом $\td{\pi}^{n}_{k}\circ\td{\pi}^{k}_{i}=\td{\pi}^{n}_{i}$,
  $n,k,i\in\mm{N}, (n>k>i),$ то мы имеем проективную последовательность
   топологических пространств
  $
  \{
  \td{\fr{U}}_n,\td{\pi}^{n}_{k}\},$
  и можно определить обратный предел \cite[III.1.5]{phillvv},
  \cite[2.5.1]{en}, в наших обозначениях его можно записать
   $\lim_{\leftarrow} \{\td{\fr{U}}_n,\td{\pi}^{n}_{k}\} \rav\td\Delta(\td{\fr{U}})\subset \td\Pi.$
   Как показано выше, $\td\Delta$ --- гомеоморфизм, поэтому
     $\td{\fr{U}}$ гомеоморфно $\td\Delta(\td{\fr{U}})$.
  Теперь в силу теоремы Куроша \cite[III.1.13]{phillvv}
   обратный предел $\td\Delta(\td{\fr{U}})$ компактов $\td{\fr{U}}_n$ --- компакт, но тогда и само $\td{\fr{U}}$ также
   компакт. Более того, поскольку $\td{\fr{U}}_n$ метризуемы, то аналогично из
  \cite[4.2.5]{en} $\td{\fr{U}}$ также и метризуемо.

  Повторив рассуждения, но уже без  $\td{\ }$,
  или просто сославшись на  \cite[3.4.11]{en} (пространство функций в  компактно-открытой топологии
  как предел обратной последовательности) и \cite[2.5.6]{en} (сужение спектра по замкнутым множествам),
  имеем   $\displaystyle{\fr{U}}\cong\lim_{\leftarrow} \{
  {\fr{U}}_n,\pi^{n}_{k}\}\rav \Delta({\fr{U}})\subset \Pi.$

   Введем отображение ${e}_n:\fr{U}_n\to \td{\fr{U}}_n$ правилом
   ${e}_n(u)(t)\rav(\td\delta\circ u)(t)=\td{\delta}_{u(t)}$ для всех
   $n\in\mm{N},t\in[0,n],u\in\fr{U}_n.$ Поскольку для всех $n,k\in\mm{N},
  n>k$  выполнено
  ${e}_k \circ {\pi}^n_k = {e}_n,$
  то семейство отображений ${e}_n$ также можно превратить в
  проективную систему $\{{e}_n,{\pi}^n_k\}$ и поднять до отображения между пределами систем
  $\fr{U}_n$ и
  $\td{\fr{U}}_n$, получим  обратный предел $e_\Delta:\Delta(\fr{U})\to\td\Delta(\td{\fr{U}})$,
  при этом,  из $e_n\circ\pi_n=\td\pi_n\circ\td\delta$ имеем
   $e_\Delta\circ\Delta=\td\Delta\circ\td{\delta},$ а из
  $\td{\fr{U}}_n=cl e_n(\fr{U}_n)$ (\cite{va}) имеем
  $\td\Delta(\td{\fr{U}})=cl e_\Delta\big(\Delta(\fr{U})\big)=cl (\td\Delta\circ\td{\delta})(\fr{U});$
  теперь, в силу непрерывности $\td\Delta^{-1}$ получаем $\td{\fr{U}}= cl \td{\delta}(\fr{U}),$
  и первый пункт доказан.

  Пусть дано некоторое неупреждающее отображение $\phi:\fr{U}\to
   C_{loc}(\mm{T},E);$
  вследствие неупреждаемости $\phi$ корректно ввести отображение
  $\phi_n:\fr{U}_n\to C([0,n],E)$
  правилом $\phi_n\big(\pi_n(u)\big)=\phi(u)$ для любого $u\in\fr{U}.$
  Поскольку для всякого $n\in\mm{N}$ $e_n(\td{\fr{U}})$ всюду плотно в $\td{\fr{U}}_n,$
  то  существует  единственное непрерывное
  отображение
  $\td\phi_n:\td{\fr{U}}_n\to C([0,n],E)$ со свойством
 $\phi_n=\td\phi_n\circ
  {e}_n$.
  Но тогда если
  определить $\td\phi$ правилом $\td\phi(u)|_{[n-1,n\ra}\rav(\td\phi_n\circ \td\pi_n)_{[n-1,n\ra}$
  для всех $n\in\mm{N},$ $u\in\fr{U}$,
  то для всех $n\in\mm{N}$
  $\td\phi|_{[0,n]}=\td\phi_n\circ \td\pi_n,$ в частности $\td\phi$
  непрерывно (топология компактно-открыта), а кроме того
  $(\td\phi\circ\td{\delta})(u)|_{[0,n]}=(\td\phi_n\circ \td\pi_n\circ\td{\delta})(u)=
  (\td\phi_n\circ e_n\circ \pi_n)(u)=(\phi_n\circ \pi_n)(u)=\phi(u)|_{[0,n]}.$
  В силу произвольности $n\in\mm{N}$ существование нужного
  отображения $\td\phi:\td{\fr{U}}\to C(\mm{T},E)$ показано.
   Осталось заметить, что если бы было отличное от $\td\phi$
   отображение с тем же свойством, то  оно бы отличалось на некотором открытом
   в силу \cite[2.4.17]{phillvv}
    множестве, в частности на некоторых элементах из $\td{\delta}(\fr{U}),$
   чего быть не может. Второй пункт доказан.

   Отображение $\td\ph$ непрерывно, например, в силу \cite[Theorem
   3.5.6]{tovst1}. Аналогично будет непрерывно зависеть от $\eta\in\td{\fr{U}}$
   решение  $\phi[\eta]=(x,y)[\eta]$ задачи Коши
   $\dot{x}=\int_{U(t)}f(t,x(t),u)\eta(du)(t),
   \dot{y}=\int_{U(t)}g(t,x(t),u)\eta(du)(t), c(0)=0,x(0)=0$
   как элемент $C_{loc}(\mm{T},\mm{X}\times\mm{R}),$
   в частности $y[\eta]\in C_{loc}(\mm{T},\mm{R}),$
   но $|\td{J}(\eta)-y[\eta](n)|<\omega(n)$ для всех $n\in\mm{N}$ в силу условия $\bf{(g)}.$
   Осталось заметить, что $y[\eta](n)$ зависит лишь от $\eta|_{[0,n]}$
   и зависит непрерывно,
   $\omega(n)\downarrow 0,$ а кроме того $|\td{J}(\eta)|$ ограничено
   числом $\omega(0).$

    Последний пункт следует из предыдущих.
\bo

 Итак, обобщенная задача
 максимизации
 функционала $\td{J}(\eta)$ имеет решение, более того
 для всякого такого решения найдется сходящаяся к ней
 (в топологии $\td{\fr{U}}$)
 максимизирующая \rref{opt} последовательность
 управлений из $\fr{U}.$

 \section{Соотношения принципа максимума}
\setcounter{equation}{0}

 Определим функцию Гамильтона-Понтрягина
  $\ct{H}:\mm{X}\times Gr {U}\times\mm{T}\times\mm{X}\mapsto\mm{R}$
  правилом:
\begin{equation}
   \label{Ham}
   \ct{H}(x,t,u,\lambda,\psi)\rav\psi f\big(t,x,u\big)+\lambda g\big(t,x,u\big).
\end{equation}
 Введём соотношения
\begin{equation}
   \label{sys_x}
       \dot{x}(t)= f\big(t,x(t),u(t)\big);
   \end{equation}
\begin{equation}
   \label{sys_psi}
       \dot{\psi}(t)\in-{\partial_x \ct{H}\big(x(t),t,u(t),\lambda,\psi(t)\big)};
   \end{equation}
\begin{equation}
   \label{maxH}
        \ct{H}\big(x(t),t,u(t),\lambda,\psi(t)\big)=\sup_{p\in
        U(t)}\ct{H}\big(x(t),t,p,\lambda,\psi(t)\big);
   \end{equation}
\begin{equation}
   \label{dob}
   x(0)=0,\ \
       \ \ ||\psi(0)||_\mm{X}+\lambda=1.
   \end{equation}
 Легко видеть, что
 для всех $u\in\fr{U}$ при любых начальных
 условиях система \rref{sys_x},\rref{sys_psi} имеет локальное
 решение,
 и всякое решение этих соотношений продолжимо на всё $\mm{T}.$

  Обозначим через $\fr{Y}$ семейство всевозможных решений
   $(x,u,\lambda,\psi)\in C_{loc}(\mm{T},\mm{X})\times\fr{U}\times[0,1]\times C_{loc}(\mm{T},\mm{X})$
    системы \rref{sys_x},\rref{sys_psi},\rref{dob} на $\mm{T}$,
     через $\fr{Z}$ --- семейство тех решений из
   $\fr{Y}$, для которых  почти всюду на $\mm{T}$ выполнено также \rref{maxH}.

\medskip

   Введём такие соотношения для  обобщенных
   управлений, а именно в условиях \rref{dob} рассмотрим
\begin{equation}
   \label{sys_x_}
       \dot{x}(t)= \int_{U(t)} f(t,x(t),u) \eta(t)(du);
   \end{equation}
\begin{equation}
   \label{sys_psi_}
       \dot{\psi}(t)\in-\int_{U(t)} {\partial_x \ct{H}\big(x(t),t,u,\lambda,\psi(t)\big)}\eta(t)(du);
   \end{equation}
\begin{equation}
   \label{maxH_}
        \int_{U(t)}\ct{H}\big(x(t),t,u,\lambda,\psi(t)\big)\eta(t)(du)=\sup_{p\in
        U(t)}\ct{H}\big(x(t),t,p,\lambda,\psi(t)\big).
   \end{equation}
 Аналогично, для всех $\eta\in\td{\fr{U}}$ при любых начальных
 условиях система \rref{sys_x_},\rref{sys_psi_} имеет локальное
 решение, и
 всякое такое решение продолжимо на всё $\mm{T}.$

   Обозначим также через $\td{\fr{Y}}$ семейство всевозможных четверок
   $(x,\eta,\lambda,\psi)\in C_{loc}(\mm{T},\mm{X})\times\td{\fr{U}}\times[0,1]\times C_{loc}(\mm{T},\mm{X})$
   решений системы
   \rref{dob}--\rref{sys_psi_}.
    Введем также
    $\td{\fr{Z}}$ --- семейство таких четверок
   $(x,\eta,\lambda,\psi)\in\td{\fr{Y}},$
   для которых  почти всюду на $\mm{T}$ выполнено также и \rref{maxH_}.

   Заметим, что    для всякого $\eta\in\td{\fr{U}}$   семейство всевозможных решений
   $(x,\eta,\lambda,\psi)\in\td{\fr{Y}}$
    системы \rref{dob}--\rref{sys_psi_}
    на $\mm{T}$
    при заданном управлении $\eta,$
   компактно в силу \cite[Theorem 3.4.2]{tovst1}.
   Более того, это множество, как многозначное зависящее от $\eta$ отображение,
   полунепрерывно сверху зависит от $\eta.$
   Действительно правая часть \rref{sys_x_},\rref{sys_psi_} выпукла, интегрально ограничена,
   полунепрерывно сверху зависит от $\eta,$ при всяком фиксированном $x,\psi$ измерима,
   следовательно имеет   измеримый селектор (\cite[Lemm
   2.3.11]{tovst1});
   более того,  все локальные решения \rref{sys_x_},\rref{sys_psi_}
   продолжимы на всё $\mm{T}$. Поскольку  выполнены все условия \cite[Theorem 3.5.6]{tovst1},
   то показана
   полунепрерывность сверху этого отображения.
   Теперь компактами являются и $\td{\fr{Y}}$, и $\td{\fr{Z}}$, как
   графики этого отображения на компактной подобласти области определения.

\medskip

  Покажем сейчас, следуя идеологии \cite[\S 9]{kr_as}, принцип максимума
   в нашей обобщенной задаче.


 Зафиксируем некоторое $\eta^0\in\td{\fr{U}}$, удовлетворяющее \rref{supmax}.
 Введем для сокращения записи $x^0\rav\td{\ph}[\eta^0].$
 Обозначим через $\td{\fr{Y}}^0$ семейство таких троек $(x^0,\lambda,\psi),$
 что $(x^0,\eta^0,\lambda,\psi)\in\td{\fr{Y}}.$
 Отметим, что как показано выше, это множество --- компакт.

   Из условия $(\bf{u})$ отображение $U$ аппроксимируется счетным числом своих
   измеримых селекторов (представление Кастена)
   \cite[Theorem 12.1]{select}. Пусть последовательность
   $(v_i)_{i\in\mm{N}}\in\fr{U}^\mm{N}$ является представлением
   Кастена для отображения $U.$

 Исключительно ради удобства обозначений
 вложим конечномерное пространство $\mm{U}$ как линейное
 подпространство в также конечномерное пространство большей
 размерности  $\mm{U}'.$ Выберем какую-нибудь точку
  $v_0\in \mm{U}'\setminus\mm{U}.$
 Доопределим для всех $(t,x)\in\mm{T}\times\mm{X}$
 $$f(t,x,v_0)\rav\int_{U(t)} f(t,x,u)\eta^0(du),\   g(t,x,v_0)\rav\int_{U(t)} g(t,x,u)\eta^0(du).$$
 Теперь можно доопределить  $\ct{H}$ при $u=v_0$ при помощи правила \rref{Ham}.

   Зафиксируем $n\in\mm{N}$,
   введем симплекс
   $\displaystyle\triangle_n\rav
  \Big\{\bar{\alpha}\rav\!(\alpha_0,\dots,\alpha_{n})\in[0,1]^{n+1}\,\Big|\,
  \sum_{i=0}^{n}\!\alpha_i=1\Big\}.$
    Пусть $\fr{L}_n$ --- множество всевозможных измеримых отображений из $\mm{T}$ в $\triangle_n.$
  Введем также  $\bar{\alpha}^0\in\fr{L}_n$ правилом $\bar{\alpha}^0(t)\rav(1,0,\dots,0)$
  для всех $t\in\mm{T}.$

   Определим, наконец, вспомогательную линейную систему
 $$   \dot{x}=\sum_{i=0}^{n}\alpha_i(t) f\big(t,x,v_i(t)\big),\ x(0)=0,\ t\in\mm{T},\
    x\in\mm{X},\ \bar{\alpha}\in\fr{L}_n;$$
   каждое $\bar{\alpha}\in\fr{L}_n$ порождает траекторию
   $\bar{\ph}^n[\bar{\alpha}]\in C(\mm{T},\mm{X})$
   (причём
    $\bar{\ph}^n[\bar{\alpha}^0]=x^0$).
   Введём  на $\fr{L}_n$
   функционал
$$   \bar{J}_n(\bar\alpha)\rav\int_{[0,n]} \sum_{i=0}^{n}\alpha_i(t)
   g\big(t,\bar{\ph}^n[\bar{\alpha}](t),v_i(t)\big)dt.$$

   Покажем с помощью принципа Беллмана, что комбинация $\bar{\alpha}^0$ доставляет максимум
   этого функционала среди всех комбинаций  $\bar{\alpha}\in\fr{L}_n,$
   удовлетворяющих краевому условию $\bar{\ph}^n[\bar{\alpha}](n)=x^0(n).$
   Действительно, краевое условие  на $\bar{\alpha}^0$ выполнено по построению, пусть
   есть $\bar{\beta}\in\fr{L}_n$, для которого
    $\bar{J}_n(\bar\beta)>\bar{J}_n(\bar\alpha^0)$ и
     $\bar{\ph}^n[\bar{\beta}](n)=\td{\ph}[\eta^0](n).$
    Тогда введём $\eta^{00}\in\td{\fr{U}}$ равенствами
        $\eta^{00}|_{[n,\infty\ra}=\eta^0|_{[n,\infty\ra}$ и
        $\eta^{00}|_{[0,n\ra}=\sum_{k=0}^n\beta_k(\td\delta\circ v_k)|_{[0,n\ra}$
    Отсюда
    $\td{\ph}[\eta^{00}]|_{[n,\infty\ra}=\td{\ph}[\eta^{0}]|_{[n,\infty\ra},$
    разбивая в определении $\td{J}$ промежуток интегрирования на
    промежутки $[0,n\ra$, $[n,\infty\ra$ имеем
    $\td{J}(\eta^{00})-\td{J}(\eta^{0})=\bar{J}_n(\bar\beta)-\bar{J}_n(\bar\alpha^0)>0,$
    что противоречит оптимальности $\eta^{0}$ в \rref{opt}. Итак, $\bar{\alpha}^0$
    оптимальна.

   Поскольку  пара $(x^0,\bar\alpha^0)$
   оптимальна в этой задаче, то
  из \cite[теорема 5.2.1]{clarke}
 при некоторых $\lambda^n\in\mm{T},\psi^n\in C_{loc}([0,n],\mm{X})$
 для $t\in[0,n]$  выполнены аналоги \rref{sys_x_},
 \rref{dob} и соотношения
$$       \dot{\psi}^n(t)\in- \sum_{i=0}^{n}\alpha^0_i(t)
       {\partial_x
       \ct{H}\big(x^0(t),t,v_0(t),\lambda^n,\psi^n(t)\big)},$$
        $$\sum_{i=0}^{n}\alpha^0_i(t)\ct{H}\big(x^0(t),t,v_i(t),\lambda^n,\psi^n(t)\big)=
        \sup_{\bar{\beta}\in
        \fr{L}_n}\sum_{i=0}^{n}\beta_i(t)\ct{H}\big(x^0(t),t,v_i(t),\lambda^n,\psi^n(t)\big),$$
        теперь, воспользовавшись определениями $\bar\alpha^0,$  $f(t,x,v_0),g(t,x,v_0),$
        имеем  кроме  \rref{sys_x_},  \rref{dob}  для почти всех $t\in[0,n]$ также
   \begin{equation}
   \label{sys_psiv}
          \dot{\psi}^n(t)\in- \int_{U(t)}
       {\partial_x
       \ct{H}\big(x^0(t),t,u,\lambda^n,\psi^n(t)\big)}\eta^0(t)(du);
   \end{equation}
   \begin{equation}
   \label{sys_maxv}
        \int_{U(t)}\ct{H}\big(x^0(t),t,u,\lambda^n,\psi^n(t)\big)\eta^0(t)(du)\geq
           \max_{i\in\overline{1,n}}
        \ct{H}\big(x^0(t),t,v_i(t),\lambda^n,\psi^n(t)\big).
   \end{equation}
 Продолжим $\psi^n$ с $[0,n]$ на $\mm{T}$
 как решение \rref{sys_psiv} произвольным образом.
 В силу \cite[Теорема 2.7.2]{clarke} при перемене
 $\partial$ и $\int$ множество в \rref{sys_psiv} лишь увеличится, то есть показано
 \rref{sys_psi_}, в частности
 $(x^0,\lambda^n,\psi^n)\in\td{\fr{Y}}^0$
 для всех $n\in\mm{N}$.

 Поскольку $\td{\fr{Y}}^0$ --- компакт,
  то, перейдя к подпоследовательности, можно считать, что
  $(x^0,\lambda^n,\psi^n)_{n\in\mm{N}}$ сходится к некоторому
  $(x^0,\lambda^\infty,\psi^\infty)\in\td{\fr{Y}}^0.$
   С другой стороны, рассмотрим произвольное $T\in\mm{T}$, при
    всяком $n\in\mm{N},n>T$ для почти всех $t\in [0,T]$ выполнено
   \rref{sys_maxv}. Перейдя к пределу при $n\to\infty$,
   имеем  для почти всех
        $t\in\mm{T}$ (в силу произвольности $T$)
$$
     \int_{U(t)}\ct{H}\big(x^0(t),t,u,\lambda^\infty,\psi^\infty(t)\big)\eta(t)(du)\geq
           \sup_{i\in\mm{N}}
        \ct{H}\big(x^0(t),t,v_i(t),\lambda^\infty,\psi^\infty(t)\big).$$
        Но по выбору $(v_i)_{i\in\mm{N}}$ множество
         $\{v_i(t)\,|\,i\in\mm{N}\}$ вcюду плотно в $U(t)$
         для почти всех $t\in\mm{T}$, в частности правая часть
         неравенства
         равна       $\sup_{u\in U(t)}
        \ct{H}\big(x^0(t),t,v_i(t),\lambda^\infty,\psi^\infty(t)\big),$
        которая заведомо не меньше левой.
        Следовательно          для почти всех $t\in\mm{T}$ имеет место равенство, то есть
         \rref{maxH_}. Таким образом,
         $(x^0,\eta^0,\lambda^\infty,\psi^\infty)\in\td{\fr{Z}}$.
   Итак, доказано

  \begin{prp}
  \label{trud}
    В условиях $\bf{(f),(g),(u)}$ существует оптимальное обобщенное
    управление  $\eta^0\in\td{\fr{U}}$
    и соответствующая ему  траектория $x^0=\td\ph[\eta^0]$, более того
    для всякой такой пары $(\eta^0,x^0)$
    при некоторых $\lambda\in\mm{T},\psi\in C_{loc}(\mm{T},\mm{X})$
    четверка
    $(x^0,\eta^0,\lambda,\psi)\in\td{\fr{Z}}$
     удовлетворяет
     соотношениям принципа максимума
      \rref{dob}--\rref{maxH_}.
\end{prp}

 Для упрощения формулировок всюду далее мы будем предполагать, что
 оптимальное управление $u^0$ реализовалось среди элементов
 $\fr{U}$, соответствующую $u^0$ траекторию обозначим через $x^0$
 (при необходимости можно всегда восстановить эти формулировки
 объявив обобщенное оптимальное некоторым допустимым
 $v^0\in\mm{U}'\setminus\mm{U}$).
  Отметим, что очень общие условия существования оптимального решения среди
 допустимых управлений имеются в \cite{bald}.



\section{Условия трансверсальности}
\setcounter{equation}{0}

  Cоотношения
  \rref{sys_x}--\rref{dob}
  не содержат условия на
  правом конце. Есть несколько вариантов таких дополнительных
  условий (подробнее см.\cite[\S 1.6]{kr_as}),
  в данной работе исследуется прежде всего вариации условия
       $\displaystyle {\lim}_{t\to\infty} ||\psi(t)||_\mm{X}=0.$
 Сформулируем предположения в терминах  устойчивости $\psi$, при которых такое условие будет
 действительно необходимым.

   Через $(Fin)(u^0)$ будем обозначать семейство тех $\eta\in\td{\fr{U}},$
   для которых
   $\eta|_{[T,\infty \rangle}=(\td\delta\circ u^0)|_{[T,\infty \rangle}$
   при некотором  $T\in\mm{T}$.

  Пусть $w:\mm{T}\times \mm{U}\to \mm{T}$ отображение
  Каратеодори, обладающее измеримая мажорантой.
  Для любых $\tau\in\mm{T}$,
  $\eta\in\td{\fr{U}}$ введём
  $\displaystyle \ct{L}_w[\eta](\tau)\rav\int_{[0,\tau]}\int_{U(t)} w(t,u)\eta(t)(du)dt.$
    Предположим, что
   $\ct{L}_w[\td\delta\circ u^0]\equiv 0$ и если
   для некоторых $\eta\in D(w)$, $\tau\in\mm{T}$
   имеет место $\ct{L}_w[\eta](\tau)=0$, то  $\eta$ почти всюду на $[0,\tau]$
   равна $\td{\delta}\circ u^0$. Множество таких $w$,
      обозначим через $(Null)(u^0)$.
   Определим
   окрестность
   $
   O_{w}(u^0;\epsi)\rav\{\eta\in (Fin)(u^0)\,|\, \forall
   t\in\mm{T}\
    \ct{L}_w[\eta](t)<\epsi\}$ для всех $w\in(Null)(u^0),\epsi\in\mm{R}_{>0}.$

 Условие $\bf{(s)}$: На  оптимальной для задачи \rref{opt} паре $(x^0,u^0)$
 существует такой вес $w\in(Null)(u^0)$, что
 для всякого решения  $(x^0,u^0,\lambda^0,\psi^0)\in\widetilde{\fr{Z}}$
 множитель Лагранжа $\psi^0$  устойчив к малым для $\ct{L}_w$ возмущениям системы
 \rref{sys_x},\rref{sys_psi},\ 
  то есть
  для всякого $\epsi\in\mm{R}_{> 0}$
  найдутся такие  число $\delta\in\mm{R}_{> 0}$  и окрестность $\Upsilon\subset\td{\fr{Y}}$
  решения $(x^0,\td\delta\circ u^0,\lambda^0,\psi^0)$ (как элемента
  $C_{loc}(\mm{T},\mm{X})\times \td{\fr{U}}\times\mm{T}\times C_{loc}(\mm{T},\mm{X})$),
  что
   для всех  $(x,\eta,\lambda,\psi)\in\Upsilon$ из
   $||\ct{L}_w[\eta]||_C<\delta$
   следует $||\psi-\psi^0||_C<\epsi$.

\begin{prp}
\label{3}
    В условиях $\bf{(f),(g),(u)}$
    для всякой оптимальной для задачи \rref{opt} пары
    $(x^0,u^0)\in C_{loc}(\mm{T},\mm{X})\times\fr{U}$,
    удовлетворяющей условию
    $\bf{(s)}$, для всякой неограниченно возрастающей
    последовательности
  моментов времени $(\tau_n)_{n\in\mm{N}}\in \mm{T}^\mm{N}$
    найдется такое решение
    $(x^0,u^0,\lambda^\infty,\psi^\infty)\in\fr{Z}$
    всех соотношений принципа максимума
    \rref{sys_x}--\rref{dob},
    что   выполнено также
    условие трансверсальности
    \begin{equation}
   \label{partlim}
       \underline{\lim}_{n\to\infty}||\psi^\infty(\tau_n)||_{\mm{X}}=0.
   \end{equation}
\end{prp}

\doc
 Зафиксируем некоторую неограниченную, монотонно возрастающую
    последовательность
  моментов времени $(\tau_n)_{n\in\mm{N}}\in \mm{T}^\mm{N}$.
  Возьмем также произвольную сходящуюся к нулю последовательность
  $(\gamma_n)_{n\in\mm{N}}\in \mm{T}^\mm{N}$ со свойством:
  $\omega(\tau_n)/\gamma_n\to 0$.
  (Например, подойдет $\gamma_n\rav \sqrt{\omega(\tau_n)},$
  где функция $\omega$ взята из условия ${\bf{(g)}}$).

  Рассмотрим для каждого $n\in\mm{N}$ задачу
\begin{equation}
   \label{z_1}
   J_n(\eta)\rav\int_{[0,\tau_n\rangle}\int_{U(t)}
  g(t,\td{\ph}[\eta](t),u)\eta(t)(du)dt-\gamma_n \ct{L}_w[\eta](\tau_n)
  \to\max.
   \end{equation}
  Функционал здесь ограничен сверху числом $J(u^0)+\omega(\tau_n),$
  следовательно имеет супремум, но каждое слагаемое
  непрерывно зависит от $\eta,$ пробегающего компакт $\td{\fr{U}}$,
   следовательно эта задача имеет в $\td{\fr{U}}$ оптимальное
  решение; какое-нибудь из них обозначим через $(x^n,\eta^n)$.
  Как при доказательстве предложения \ref{trud} можно записать
  $\eta_n$ как фиктивное допустимое управление, зафиксировать разложение Кастена и
  рассмотреть линейную вспомогательную задачу  для ${J}_n$ при  $n\in\mm{N}.$
  Тогда,
  переходя к пределу при $n\to\infty$, получим, что
  при каких-то $(\lambda^n,\psi^n)\in\mm{T}\times C([0,n],\mm{X})$ на промежутке $[0,\tau_n\rangle$
  имеет место
  \rref{dob}--\rref{sys_psi_}, а
  для отображения
  $\displaystyle \ct{H}_{\tau_n}(x,t,u,\lambda,\psi)\rav\ct{H}(x,t,u,\lambda,\psi)-
  \gamma_n\int_{U(t)}w(t,u)\eta(t)(du)$
  выполнено для почти всех $t\in\mm{T}$
\begin{equation}
   \label{sys_max_}
   \int_{U(t)}\ct{H}_{\tau_n}\big(x(t),t,u,\lambda,\psi(t)\big)
   \eta(t)(du)=
   \sup_{p\in
        U(t)}\ct{H}_{\tau_n}\big(x(t),t,p,\lambda,\psi(t)\big).
   \end{equation}
  Более того, $(x^n,\eta^n,\lambda^n,\psi^n)$ удовлетворяет
  условию трансверсальности на свободном конце: $\psi^n(\tau_n)=0.$
  На $[\tau_n,\infty\rangle$
   всю четверку $(x^n,\eta^n,\lambda^n,\psi^n)$
  продолжим при помощи управления $u^0|_{[\tau_n,\infty\rangle}.$
  Тогда $\eta^n\in (Fin)(u^0).$
  Обозначим через $\fr{Z}^n$ множество четверок
  $(x,u,\lambda,\psi)$,
  удовлетворяющих на $t\in\mm{T}$
  соотношениям
  \rref{dob}--\rref{sys_psi_},
  соотношению   \rref{sys_max_} почти всюду
   на $[0,\tau_n\rangle$, и со свойством $(\td{\delta}\circ u^0)|_{[\tau_n,\infty\rangle}=
   \eta^n|_{[\tau_n,\infty\rangle}.$ Теперь $(x^n,\eta^n,\lambda^n,\psi^n)\in\fr{Z}^n$
   для всякого $n\in\mm{N}.$

 Заметим, что все $\fr{Z}^n$ замкнуты, а поскольку содержатся в компакте $\td{\fr{Y}}$,
  то и компактны; теперь последовательность
  $(x^n,\eta^n,\lambda^n,\psi^n)_{n\in\mm{N}}$ имеет предельную точку $(x^\infty,\eta^\infty,\lambda^\infty,\psi^\infty)\in\td{\fr{Y}}.$
  Перейдя при необходимости к подпоследовательности будем считать эту точку
  пределом всей
  последовательности.

 Для всякого фиксированного $x$ множество тех $u\in U(t),$ что реализуют
 в \rref{sys_max_}  максимум, имеет селектор в силу
 \cite[Theorem 3.7]{select}. Тогда из \cite[Lemm 2.3.11]{tovst1}
 он имеется и при подстановке  в $\ct{H}$ произвольной непрерывной функции
 $x$.
  Поскольку, кроме того, соотношение \rref{sys_max_} также полунепрерывно
  сверху зависит от $x,\psi$ и параметров $\gamma,\lambda$, а
   все соотношения интегрально ограничены на ограниченных множествах, то в силу
  \cite[Theorem 3.5.6]{tovst1}
  на каждом конечном промежутке
  для  пучков решений
   \rref{sys_x},\rref{sys_psi}, удовлетворяющих \rref{sys_max_},
    имеет место
 полунепрерывность сверху  уже
   по  $\gamma,\lambda$.\ 
   В частности, для  $\gamma\to 0$, мы получаем, что
     верхний предел компактов $\fr{Z}^n$ вложен в $\td{\fr{Z}}.$
   Отсюда
   $(x^\infty,\eta^\infty,\lambda^\infty,\psi^\infty)\in\td{\fr{Z}}.$

   Далее, в силу оптимальности  $u^n$  и $u^0$ в своих задачах, а также условия
   $\bf{(g)}$,
   имеет место
$$   J(u^n)+\omega(\tau_n)\geq J_n(u^n)\geq
   J_n(u^0)\geq J(u^0)-\omega(\tau_n)\geq
   J(u^n)-\omega(\tau_n),$$
   тогда
   $\displaystyle
      \gamma_n \ct{L}_w[\eta^n](\tau_n)=\gamma_n
     \int_{[0,\tau_n\rangle}\int_{U(t)}
  w(t,u)\eta^n(t)(du)  dt\leq 2 \omega(\tau_n),$
    и в силу
    $(\td{\delta}\circ u^0)|_{[\tau_n,\infty\rangle}=
   \eta^n|_{[\tau_n,\infty\rangle}$
\begin{equation}
   \label{to_w}
   \ct{L}_w[\eta^n](\tau)\leq 2 \omega(\tau_n)/\gamma_n
   \ \forall \tau\in\mm{T}.
\end{equation}
     Переходя для всякого $\tau\in\mm{T}$ к пределу при $n\to\infty$ имеем
   $\ct{L}_w[\eta^\infty]\leq 0$, то есть
   $\ct{L}_w[\eta^\infty](\tau)=0$ для всех $\tau\in\mm{T}.$
    Поскольку  $w\in(Null)(u^0)$, то $\eta^\infty=\td{\delta}\circ u^0$
    почти всюду на $\mm{T}$, отсюда
   $x^\infty=x^0$ и $(x^0,u^0,\lambda^\infty,\psi^\infty)\in{\fr{Z}}.$
   Более того
   из \rref{to_w}
 $||\ct{L}_w[\eta^n]||_C\to 0$.

  Выберем некоторое $\epsi\in\mm{R}_{> 0}$
  и из условия $\bf{(s)}$ возьмем $\Upsilon\subset\td{\fr{Y}}$,
  $\delta\in\mm{R}_{> 0}$, тогда найдется $N\in\mm{N}$, что
  при $n\in\mm{N}, n>N$
  выполнено
  $(x^n,\eta^n,\lambda^n,\psi^n)\in\Upsilon$,
   $||\ct{L}_w[\eta^n]||_C<\delta;$  теперь из
   условия $\bf{(s)}$ следует и
      $||\psi^n(\tau_n)-\psi^0(\tau_n)||_\mm{X}<\epsi$; но
      $\psi^n(\tau_n)=0$, откуда
   $||\psi^0(\tau_n)||_\mm{X}<\epsi$ для всех $n\in\mm{N}, n>N$.
  В силу произвольности $\epsi\in\mm{R}_{> 0}$
  показано \rref{partlim}.
\bo

Одно из самых общих условий на \rref{partlim} показано в \cite{norv}.
 Для  задачи управления без фазовых ограничений результат \cite[Theorem 6.1]{norv}
  следует из предложения 1 и
 \cite[Lemm 3.1]{norv} .

 \medskip

    \subsubsection*{Условие ${\bf{(s)}}$ как следствие частичной устойчивости}

   Цель данного раздела --- подобрать такой вес $w^0\in(Null)(u^0)$, чтобы
    условие $\bf{(s)}$ следовало бы из некоторого
    варианта (неасимптотической) устойчивости по Ляпунову компоненты $\psi$.

  Всюду далее предполагаем, что выполнено

 Условие ${\bf{(\partial)}}:$
      существуют производные
     $\frac{\partial f(t,x,u)}{\partial x},$
     $\frac{\partial g(t,x,u)}{\partial x},$
     являющиеся на $\mm{T}\times\mm{X}\times\mm{U}$
     локально липшицевыми отображениями Каратеодори,
      имеющими на каждом компакте
  суммируемую мажоранту.

  Введем для удобства ссылок
  \begin{equation}
   \label{sys_xx}
       \dot{x}(t)= f\big(t,x(t),u^0(t)\big);
   \end{equation}
\begin{equation}
   \label{sys_psii}
       \dot{\psi}(t)=-{\partial_x \ct{H}\big(x(t),t,u^0(t),\lambda,\psi(t)\big)};
   \end{equation}
\begin{equation}
\label{sys_lambda}
  \dot{\lambda}=0.
\end{equation}

 \medskip

Пусть  дано конечномерное евклидово пространство
  $E$, задано ограниченное на всяком компактном множестве мультиотображение
  $G:\mm{T}\to (comp)(E)$,
  график $Gr\, G$ которого\ ---\  замкнутое борелевское множество в  $\mm{T}\times E.$
  Пусть также дана локально липшицевая  функция Каратеодори
  $a:\mm{T}\times{E}\times \mm{U}\mapsto E$,
   имеющая измеримую мажоранту и удовлетворяющая условию
   продолжимости.
 Рассмотрим систему
\begin{equation}
   \label{a}
   \dot{y}=a(t,y(t),u),\
   t\in \mm{T}, u\in U(t).
\end{equation}

    Отметим, что также как от \rref{sys} переходили к \rref{sys_x_},
    можно рассмотреть задачу Коши
    $$\dot{y}=\int_{U(t)}a(t,y(t),u)\eta(t)(du),\ y(0)=\xi\in E, t\in \mm{T} $$
    для всех $\eta\in\td{\fr{U}},\xi\in E.$
    Эта задача вновь будет иметь единственное решение $y^\eta_\xi\in C_{loc}(\mm{T},E).$

    Для всякой позиции $(t^*,y^*)\in\mm{T}\times{E}$ существует
    решение $Y[t^*,y^*]$ задачи Коши
    $$\dot{y}=a\big(t,y(t),u^0(t)\big),\ y(t^*)=y^*, t\in \mm{T},$$
    оно
    продолжимо на всё $\mm{T}$ и
    непрерывно (как элемент $C_{loc}(\mm{T},\mm{X})$) зависит от $(t^*,y^*).$
    Введем
    $\beg(t^*,y^*)\rav Y[t^*,y^*](0)$
    для всякого $(t^*,y^*)\in\mm{T}\times E,$ при этом $y^{u^0}_{\beg(t^*,y^*)}=Y[t^*,y^*].$

   Следующее предложение предъявит метод подбора веса $w^0$.
 Заметим, что в отличие от предыдущих параграфов,
 при его доказательстве, а следовательно и всюду ниже, условие
 продолжимости формально нужно не только для продолжения локального
 решения  вправо до бесконечности, но и для продолжения
 локальных решений вплоть до нуля.

\begin{prp}
\label{dop} Пусть выполнено $\bf{(u)}.$
   Пусть мультиотображение
  $G:\mm{T}\to (comp)(E)$ ограничено на всяком компактном
  множестве, а его
  график $Gr\, G$ --- замкнутое борелевское множество в  $\mm{T}\times E.$
  Пусть $a:\mm{T}\times{E}\times \mm{U}\mapsto E$ --- локально липшицевая  функция Каратеодори,
   имеющая измеримую мажоранту и удовлетворяющая условию
   продолжимости.

    Тогда для всякой  $u^0\in\fr{U}$
    найдётся такое $w^0\in(Null)(u^0),$ 
         что для  произвольных $\eta\in\td{\fr{U}}, \xi\in E$,
         $T\in\mm{T}$ для $y=y^\eta_\xi$
   из $Gr\, y|_{[0,T]}\subset Gr\, G$ следует 
 $$  \big|\big|\beg\big(\tau,y(\tau)\big)-y(0)\big|\big|_E\leq \ct{L}_{w^0}[\eta](\tau)\ \ \ \forall \tau\in [0,T].$$
 \end{prp}
\doc
      Зафиксируем $n\in\mm{N}$. В силу условия продолжимости всякому $(t^*,y^*)\in Gr\,\ G|_{[0,n]}$
      можно сопоставить траекторию $Y[t^*,y^*]|_{[0,n]}$, а значит и
      число $\beg(t^*,y^*).$ В силу теоремы о непрерывной зависимости
      от начальных данных эта функция будет непрерывна, в частности,
      будет замкнут и его образ --- множество
      $$\bar{G}_n\rav\bigcup_{(t^*,y^*)\in Gr\,\, G|_{[0,n]}} Gr\,\,
      y^{u^0}_{\beg(t^*,y^*)}|_{[0,n]},$$
      а из условия продолжимости это множество ограничено, следовательно\ ---\ компактно.
      Тогда на этом множестве функция $a$ липшицева по $y$ для некоторой константы
      Липшица $L_n\rav L^a_{\bar{G}_n}\in L^1_{loc}(\mm{T},\mm{T}).$
      Примем для всех $t\in[0,n]$
       $M_n(t)=\int_{[0,t]}L_n(\tau)d\tau.$
       Отметим, что эта функция абсолютно непрерывна и монотонно не убывает.

    Рассмотрим для некоторого $n\in\mm{N}$
    для всех $t\in [n-1,n\ra, u\in \mm{U}$ число
    \begin{equation}
   \label{def_S}
     S(t,u)\rav \sup_{ y\in \bar{G}_n}
     \big|\big|a(t,y,u)-a(t,y,u^0(t)\big)\big|\big|_E.
   \end{equation}
 Заметим, что  норма внутри супремума --- отображение Каратеодори,
  $y$
  пробегает  компакт, теперь
 для всякого $u\in \mm{U}$
 из \cite[Theorem 3.7]{select} супремум  достигает максимума
 при
 подстановке
 некоторой функции $y_{max}(u)\in B([n,n-1\ra,\bar{G}_n)$,
  тогда $S(t,u)$  измерима по $t$ при всяком $u\in\mm{U}$.

 Зафиксируем $t\in [n-1,n\ra,$
 найдется функция $\omega^t\in\Omega$, для которой
\begin{equation}
  \label{1111}
   \Big|\,\big|\big|a(t,y,u')-a(t,y,u^0(t)\big)\big|\big|_E-\big|\big|a(t,y,u'')-a\big(t,y,u^0(t)\big)\big|\big|_E\,\Big|<
   \omega^t\big({\textstyle\frac{1}{||u'\!-\!u''||}}\big)
\end{equation}
 выполнено при любых $y\in \bar{G}_n,u',u''\in U(t) (u'\neq u'').$
  Рассмотрим произвольные различные $u',u''\in U(t)$, без ограничения общности
  считаем $S(t,u')< S(t,u'')$. 
  Теперь, с одной стороны
  $S(t,u')\geq
  \big|\big|a(t,y^{u''},u')-a\big(t,y^{u''},u^0(t)\big)\big|\big|_E$,
   подставляя  $y\rav y_{max}(u'')(t)$ в \rref{1111} с другой, имеем
  $0< S(t,u'')-S(t,u')\leq
  \omega^t(1/||u'-u''||),$
  то есть $S$ непрерывна на $Gr U|_{[n-1,n\ra}$ по переменной $u.$

   Таким образом функция $S: Gr U|_{[n-1,n\ra}\to\mm{T}$ является функцией
 Каратеодори. Заметим, что, пробегая $n\in\mm{N}$, мы определим
   функцию Каратеодори $S$ на всём $Gr U$.
 Более того, $S(t,u^0(t))\equiv 0$ по построению.
    Тогда корректно определить  $w^0\in (Null)(u^0)$ правилом:
    для всех $n\in\mm{N},(t,u)\in Gr\, U|_{[n-1,n\ra}$
 $$ w^0(t,u)\rav||u-u^0(t)||+e^{M_n(t)}S(t,u).$$

        Рассмотрим произвольные $n\in\mm{N}$,$\tau\in [0,n],$ $(\tau,y_1),(\tau,y_2)\in
        \bar{G}_n$. Тогда выполнено
        $Gr\, Y[\tau,y_1]|_{[0,n]},Gr\, Y[\tau,y_2]|_{[0,n]}\subset \bar{G}_n,$
        введем  на $[0,n]$ функции
        $$r(t)\rav Y[\tau,y_1](t)-Y[\tau,y_2](t),\
         V_+(t)\rav e^{M_n(t)}||r(t)||_E\ \forall t\in [0,n].$$
       Теперь из липшицевости $a$ имеем 
$  ||\dot{r}(t)||_E\geq -L_n(t)
       ||r(t)||_E,$ откуда
       $$\frac{dV^2_+(t)}{dt}=2L_n(t)V^2_+(t)
       +2e^{2M_n(t)}r(t)^T\dot{r}(t)
       \geq2L_n(t)V^2_+(t)
       -
       2L_n(t)V^2_+(t)
       =0.
       $$
       Тогда $V_+$ не убывает, и для всех $(\tau,y_1),(\tau,y_2)\in
        \bar{G}_n$ выполнено
\begin{equation}
  \label{1037_}
  ||\beg(\tau,y_1)-\beg(\tau,y_2)||_E={V_+(0)}\leq {V_+(\tau)}= e^{M_n(\tau)}||y_1-y_2||_E.
\end{equation}

     Рассмотрим теперь такие произвольные $\eta\in\td{\fr{U}}, \xi\in E$,
         $T\in\mm{T}$, что  для $y\rav y^\eta_\xi$
  выполнено $Gr\, y|_{[0,T]}\subset Gr\, G.$
        Зафиксируем произвольные $n\in\mm{N}$, $\tau_1,\tau_2\in[0,T]\cap[n-1,n\ra,$
        $\tau_1<\tau_2.$
        Примем
        $$y_1(t)\rav Y[\tau_1,y^\eta_\xi(\tau_1)](t),\ y_2(t)\rav
        y^\eta_\xi(t),\
        r(t)\rav y_1(t)-y_2(t),\ V_-(t)\rav e^{-M_n(t)}||r(t)||_E\ \forall t\in [\tau_1,\tau_2].$$
        По построению $\bar{G}_n$ выполнено
        $Gr\, y_1,Gr y_2\subset \bar{G}_n.$
        Теперь для почти всех $t\in[\tau_1,\tau_2]$
       $$\frac{dV_-^2(t)}{dt}
       =            -2L_n(t)V_-^2(t)+
       2e^{-2M_n(t)}{r}(t)\big(\dot{y}_1(t)-a(t,{y}_2(t),u^0(t))+
       a(t,{y}_2(t),u^0(t))-\dot{y}_2(t)\big)
       \leqref{def_S}$$ $$-2L_n(t)V_-^2(t)+
          2L_n(t)V_-^2(t)+2e^{-2M_n(t)}||r(t)||_E
          \int_{U(t)}S(t,u)
          \eta(t)(du)\leq $$ $$2e^{-M_n(t)}V_-(t)
          \int_{U(t)}S(t,u)
          \eta(t)(du)\leq
          2e^{-2M_n(t_1)}V_-(t)\frac{d
          \ct{L}_{w^0}[\eta](t)}{dt}.$$
Отметим, что $V_-(t_1)=0,$ функция $V_-$ неотрицательна, тогда
  $V_-^2$ не превосходит верхнего решения $\bar{V}$ задачи
  $$\frac{dV(t)}{dt}=2e^{-2M_n(t_1)}\sqrt{V(t)}\,\frac{d
          \ct{L}_{w^0}[\eta](t)}{dt}\ \ V(t_1)=0,$$
   откуда $$||r(\tau_2)||_E=e^{M_n(\tau_2)}V_-(\tau_2)\leq
   e^{M_n(\tau_2)}\sqrt{\bar{V}(\tau)}=e^{M_n(\tau_2)-2M_n(\tau_1)}
          \big(\ct{L}_{w^0}[\eta](\tau_2)-\ct{L}_{w^0}[\eta](\tau_1) \big).$$
          Теперь из
         $\beg(\tau_2,y_1(\tau_2))=\beg(\tau_2,Y[\tau_1,y^\eta_\xi(\tau_1)](\tau_2))=
         Y[\tau_1,y^\eta_\xi(\tau_1)](0)=\beg(\tau_1,y^\eta_\xi(\tau_1)$
  получаем
 \begin{equation}
  \label{1649}
     ||\beg(\tau_2,y^\eta_\xi(\tau_2))-\beg(\tau_1,y^\eta_\xi(\tau_1))||_E\leqref{1037_}
    e^{2M_n(\tau_2)-2M_n(\tau_1)}
          \big(\ct{L}_{w^0}[\eta](\tau_2)-\ct{L}_{w^0}[\eta](\tau_1) \big).
\end{equation}

    Зафиксируем произвольное $t\in [0,T].$
    Для всякого $\epsi\in\mm{R}_{>0}$ промежуток $[0,t\ra$ можно
    разбить на такие полуинтервалы вида $[\tau',\tau''\ra$
    так, чтобы всегда $[\tau',\tau''\ra\subset[n-1,n\ra$
    для некоторого $n\in\mm{N}$ и
    $M_n(\tau'')-M_n(\tau')=\int_{[\tau',\tau''\ra} L_n(t)dt<\epsi.$
   Тогда  для $\tau_1\rav\tau',\tau_2\rav\tau''$  выполнено \rref{1649}, то есть
$$||\beg(\tau'',y^\eta_\xi(\tau''))-\beg(\tau',y^\eta_\xi(\tau'))||_E\leq
    e^{2\epsi}
          \big(\ct{L}_{w^0}[\eta](\tau'')-\ct{L}_{w^0}[\eta](\tau') \big).$$
    Складывая по всем полуинтервалам, в силу $\beg(0,y^\eta_\xi(0))=\xi$
    и неравенства треугольника,
    имеем           для всех $t\in [0,T]$
    $||\beg(t,y^\eta_\xi(t))-\xi||_E\leq
    e^{2\epsi}
          \ct{L}_{w^0}[\eta](t)$, осталось перейти к пределу при $\epsi\to 0.$
%
 \bo


\medskip

 Пусть $E$ можно представить в виде $E=E_p\times E_q$
 для некоторых конечномерных евклидовых пространств $E_p, E_q$,
 обозначим проекции отображения $a$ на подпространства $E_p$ и $E_q$
 через $b$ и $c$ соответственно, теперь
 систему \rref{a} можно записать в виде
\begin{equation}
   \label{aa}
   \dot{p}=b(t,p,q,u), \dot{q}=c(t,p,q,u),\
 (p,q)(0)=\xi\in E,\ u\in U(t).
\end{equation}
  Тогда можно считать, что  $y^\eta_\xi\rav(p^\eta_\xi,q^\eta_\xi)$
   для всех $\eta\in\td{\fr{U}}, \xi\in E$.

 Пусть дано  замкнутое множество $G_0\subset E$ и $\xi\in G_0$.
 Будем говорить, что решение $y^{u^0}_{\xi}$ уравнения\ \rref{aa}
 имеет в $G_0$  для управления $u^0$ устойчивую  по Ляпунову компоненту $p$, если
  для любого $\epsi\in\mm{R}_{>0}$ существует
  такое $\delta(\epsi,y)\in\mm{R}_{>0}$, что для всех
  $\xi'\in G_0$
   из
  $||\xi'-\xi||_{E}<\delta(\epsi,y)$ следует
  $||p^{u^0}_{\xi'}(s)-p^{u^0}_\xi(s)||_E<\epsi$
    для всех
   $s\in\mm{T}$.

\begin{prp}
\label{zvez}
  Пусть выполнено $\bf{(u)}.$ Пусть $a:\mm{T}\times{E}\times \mm{U}\mapsto E$ --- локально липшицевая  функция Каратеодори,
   имеющая измеримую мажоранту и удовлетворяющая условию
   продолжимости.

   Пусть дано замкнутое множество $G_0\subset E$ и  такой компакт $K_0\in(comp)(G_0),$
    что для любого $\xi\in K_0$  решение $y^{u^0}_{\xi}$
 имеет в $G_0$  для управления $u^0$ устойчивую  по Ляпунову компоненту $p$.

  Тогда
  для любого $\epsi\in\mm{R}_{>0}$
        найдется такое $\delta\in\mm{R}_{>0}$, 
        что для всех $\xi\in K_0,
        \eta\in O_{w^0}(u^0;\delta)$ 
        если $\beg(t,y^\eta_{\xi}(t))\in G_0$ для всех $t\in \mm{T},$ то
        $||p^\eta_\xi-p^{u^0}_{\xi}||_C<\epsi.$
\end{prp}
 \doc
    Введем компакт $K_>\rav\{\xi\in G_0\,|\,\exists \xi_0\in K_0\,||\xi-\xi_0||_E\leq 1\}.$
    Сопоставим всякому $t\in\mm{T}$ множество
    $G(t)\rav\{y^\eta_\xi(t)\,|\,\xi\in K_>\}.$
    Получившееся отображение $G$ компактнозначно и непрерывно, в
    частности имеет замкнутый график.
    Теперь для мультиотображения $G$  можно из
     предложения \ref{dop} найти вес $w^0\in(Null)(u^0).$

  Определим
  $M(\xi',\xi'')\rav\sup_{t\in\mm{T}}
||p^{u^0}_{\xi'}(t)-p^{u^0}_{\xi''}(t)||_{E_p}\in\mm{T}\cup\{+\infty\}$
  для всяких $\xi',\xi''\in K_{>} \times K_{>}.$
  Для каждого $\xi\in K_{0}$ устойчивость компоненты $p$ означает конечность и непрерывность
  отображения $M$
    в точке $(\xi,\xi)\in K_{>} \times K_{>}.$

  Зафиксируем  $\epsi\in\mm{R}_{>0}$,
  выберем для каждого  $\xi\in  K_{0}$ своё $\delta(\epsi/2,y^{u^0}_\xi)\in\la 0,1/2]$,
  а значит и $\delta(\epsi/2,y^{u^0}_\xi)-$окрестность точки $(\xi,\xi)$
  (в $ K_>\times K_>$).
  Из полученного покрытия диагонали $\Delta$ множества $K_0\times K_0$ выделим конечное подпокрытие,
  оно порождает некоторую открытую окрестность $\Upsilon$ диагонали $\Delta$.
  Пусть $\delta(K_0)$ --- минимальное расстояние от $\Delta$
  до границы   окрестности $\Upsilon$.
  Теперь для  всех $\xi'\in K_>$, $\xi\in{K}_0$
  из $||\xi'-\xi||_E<\delta(K_0)$ следует $(\xi',\xi)\in\Upsilon,$ то есть
  для некоторого $\xi''\in K_0$ выполнено  $M(\xi,\xi''),M(\xi'',\xi')<\epsi/2,$
  откуда $M(\xi,\xi')<\epsi.$ Итак,
\begin{equation}
\label{1291} (||\xi'-\xi||_E<\delta(K_0))\Rightarrow
  (||p^{u^0}_{\xi'}-p^{u^0}_\xi)||_C<\epsi)\ \ \forall \xi\in{K}_0,\xi'\in K_>.
\end{equation}

  Рассмотрим произвольные $\eta\in O_{w^0}(u^0;\delta(K_0)),\xi\in K_0$
  со свойством: $\xi_1(t)\rav \beg(t,{y}^{\eta}_\xi(t))\in G_0$ для всех
  $t\in\mm{T}.$
  В силу $K_0\subset K_>=G(0)$ корректно определить пусть даже
  бесконечное
  $T_0\rav\sup \{T\in\mm{T}\,|\, Gr\, \xi_1|_{[0,t]}\subset K_>\ \forall t\in[0,T\ra\}
  \in \mm{T}\cup\{+\infty\}, $
  тогда $Gr\, {y}^{\eta}_\xi|_{[0,T_0\ra} \subset Gr\, G.$
  Теперь из предложения \ref{dop}  имеем 
\begin{equation}
\label{1364}
||{\xi}_1(t)-\xi||_E=||\beg(t,{y}^{\eta}_{\xi_1}(t))-\xi||_E\leq
\ct{L}_{w^0}[\eta](t)<\delta(K_0)\ \forall t\in [0,T_0\ra.
\end{equation}
  Для каждого $t\in [0,T_0\ra$ подставим  ${\xi}_1(t)\in K_>$ в \rref{1291},
  в силу
  ${p}^{u^0}_{\xi_1(t)}(t)={p}^{\eta}_\xi(t)$
    имеем
   $||{p}^{\eta}_\xi(t)-{p}^{u^0}_\xi(t)||_E<\epsi$ для всех $t\in[0,T_0\ra.$
   Для завершения доказательства осталось показать, что $T_0=\infty.$

   Пусть не так, тогда $G_0 \supset Gr\, \xi_1|_{\la T_0,\tau]}\not\subset Gr\, K_>$
    для любого $\tau\in \la T_0,\infty\ra$ по построению $T_0$,
   то есть $\xi_1(T_0)\in cl (G_0\setminus K_>)$ и  $||\xi_1(T_0)-\xi||_E\geq 1$
   по определению $K_>.$
   Но переходя к пределу в  \rref{1364} имеем
   $||\xi_1(T_0)-\xi||_E\leq\delta(K_0)\leq 1/2.$ Это противоречие доказывает, что
   $T_0=\infty$. \bo

\begin{cor}
\label{s3}
  В условиях $\bf{(f),(g),(u),(\partial)}$ пусть на оптимальной для задачи \rref{opt} паре
    $(x^0,u^0)\in C_{loc}(\mm{T},\mm{X})\times\fr{U}$ для управления $u^0$ всякое
    решение $(x^0,\psi^0,\lambda^0)$
 системы
     \rref{sys_xx}--\rref{sys_lambda}
     имеет устойчивую  по Ляпунову компоненту $\psi^0$
      в
 $G_0\rav\mm{X}\times\mm{X}\times[0,1]$.

    Тогда
    для всякой неограниченно возрастающей
    последовательности
  моментов времени $(\tau_n)_{n\in\mm{N}}\in \mm{T}^\mm{N}$
    найдется такое решение
    $(x^0,u^0,\lambda^\infty,\psi^\infty)\in\fr{Z}$
    всех соотношений принципа максимума
    \rref{sys_x}--\rref{dob}
     что   выполнено также
    условие \rref{partlim}.
\end{cor}
  \doc
   Достаточно в \rref{aa}  принять $E_p\rav\mm{X},$ $E_q\rav\mm{X}\times\mm{R},$
  под $p$ и $q=(q_1,q_2)$ понимать $\psi$ и $(x,\lambda),$
    в качестве $b$ взять правую часть \rref{sys_psii},
    в качестве $c$ взять правую часть системы из \rref{sys_xx} и
    \rref{sys_lambda}, принять
    $K_0\rav \mm{D}\times\{0_\mm{X}\}\times[0,1],$
    $G_0\rav\mm{X}\times\mm{X}\times[0,1].$
  Теперь предложение \ref{zvez} гарантирует выполнение условий предложения
  \ref{3}, что и требуется. \bo

    Отметим, что в следствии \ref{s3} можно принять
    $G_0\rav\{\beg(t,\psi(t),x(t),\lambda)\,|\,
         (x,\eta,\lambda,\psi)\in\td{\fr{Y}},t\in\mm{T}\}.$

\subsubsection*{Модификации условия трансверсальности \rref{partlim}}

 В некоторых случаях нет устойчивости множителя Лагранжа $\psi$,
 однако известна скорость его роста
  или есть устойчивость по некоторым компонентам векторной переменной $\psi$.
 В этом случае может помочь следующая модификация условия \rref{partlim}, условие
 типа
     \begin{equation}
   \label{partlim_}
       \underline{\lim}_{n\to\infty}
       ||\psi^0(\tau_n)A(\tau_n)||_\mm{X}=0,
   \end{equation}
 где $A$ --- измеримое по $t$ отображение из $\mm{T}$ в $\mm{L}$,
 $\mm{L}$\ ---\  линейное пространство всевозможных $m\times m$-матриц.
  Снарядим $\mm{L}$, для определенности
 операторной нормой.

 Примерами таких отображений могут быть например:
 сопоставляющее единичную матрицу отображение $A(t)\equiv1_\mm{L};$
 какой-либо "скалярный"\ множитель $A(t)\equiv r(t)1_{\mm{L}},$
 какое-либо отображение $A(t)\equiv D$ с диагональной матрицей $D$;
 к  такому виду  также сводится часто используемое условие $\psi(t)x(t)\to 0$.

  Пусть для всех $\eta\in\td{\fr{U}},\xi\in\mm{X}$ выбрано измеримое  отображение
  $A^\eta_\xi:\mm{T}\to\mm{L}$.
  Примем $A\rav A^{u^0}_{0}.$

  Условие $\bf{(s'A)}$: На  оптимальной для задачи \rref{opt} паре $(x^0,u^0)$
  найдется такой вес $w\in(Null)(u^0)$, что
 у всякого решения  $(x^0,u^0,\lambda^0,\psi^0)\in\widetilde{\fr{Z}}$
  произведение $A\psi$  устойчиво к малым для $\ct{L}_w$ возмущениям системы
 \rref{sys_x},\rref{sys_psi},\ 
  то есть
  для всякого $\epsi\in\mm{R}_{> 0}$
  найдутся такие  число $\delta\in\mm{R}_{> 0}$  и окрестность $\Upsilon\subset\td{\fr{Y}}$
  решения $(x^0,\td{\delta}\circ u^0,\lambda^0,\psi^0)$ (как элемента
  $C_{loc}(\mm{T},\mm{X})\times \td{\fr{U}}\times\mm{T}\times C_{loc}(\mm{T},\mm{X})$),
  что
   для всех  $(x,\eta,\lambda,\psi)\in\Upsilon$ из $||\ct{L}_w[\eta]||_C<\delta$
   следует $||\psi A^\eta_{\psi(0)}-\psi^0 A||_C<\epsi$.

\begin{prp}
\label{3_}
    В условиях $\bf{(f),(g),(u)}$
    для всякой оптимальной в задаче \rref{opt} пары
    $(x^0,u^0)$,
    удовлетворяющей условию
    $\bf{(s'A)}$, для всякой неограниченно возрастающей
    последовательности
  моментов времени $(\tau_n)_{n\in\mm{N}}\in \mm{T}^\mm{N}$
    найдется такое решение
    $(x^0,u^0,\lambda^0,\psi^0)\in\fr{Z}$
    всех соотношений принципа максимума
    \rref{sys_x}--\rref{dob},
    что   выполнено
    условие  \rref{partlim_}.
\end{prp}
  Доказательство отличается от доказательства предложения \ref{3}
  лишь ccылками на $\bf{(s'A)}$ вместо $\bf{(s)}$ и добавлением
  $A^\eta_{\xi},A$  как множителей
  в неравенства
  предпоследней строчки.\bo

\medskip

 Примем $E_p\rav\mm{X},$ $E_q\rav\mm{X}\times\mm{X}\times\mm{R},$
  под $q$ будем понимать $(x,\psi,\lambda),$
    в качестве $c$ возьмем правую часть \rref{sys_x}, \rref{sys_psi} и
    \rref{sys_lambda}, теперь используем имеющийся произвол выбора
    $b$.
\begin{zam}
\label{s3_}
    Пусть  выполнены $\bf{(f),(g),(u),(\partial)};$ задано локально липшицевое, имеющее измеримую
   мажоранту отображение Каратеодори
   $b:\mm{T}\times\mm{X}\times E_q\times \mm{U}\to\mm{X},$
   удовлетворяющее условию продолжимости;
   пусть также для всякого $u\in{\fr{U}}$ имеется функция Каратеодори
   $(t,\xi)\in\mm{T}\times E_q\mapsto A^{u}_\xi(t)\in\mm{L}$
    такая, что непрерывно отображение $(u,\xi)\in\fr{U}\times E_q \mapsto A^{u}_\xi(0)\in\mm{L}$
    Предположим также, что при всяком $u\in\fr{U}$ для всякого решения
    $z=(x,\lambda,\psi)$ системы
    \rref{sys_x},\rref{sys_psi},\rref{sys_lambda}
    с начальными условиями $\xi\in E_p$
    произведение
    $\psi A^{u}_{\xi}$
    является решением уравнения
\begin{equation}
 \label{1460}
 \frac{dp}{dt}= b(t,p(t),\psi(t),x(t),\lambda,u(t)).
 \end{equation}

   Пусть на оптимальной для задачи \rref{opt} паре
    $(x^0,u^0)\in C_{loc}(\mm{T},\mm{X})\times\fr{U}$  для всякого
     решения $y=(p,x,\psi,\lambda)$
 системы
    \rref{sys_xx}--\rref{sys_lambda},
    \rref{1460} при $u=u^0$,
    если только
    $y(0)\in \big\{(\psi_0A^{u^0}_\xi(0),0_\mm{X},\lambda,\psi_0)\,\big|\,\psi_0\in\mm{D}, \lambda\in [0,1]\big\},$
    то решение $y$
     имеет устойчивую  по Ляпунову компоненту $p$
     для управления $u^0$ и начальных условий из
     $\mm{X}\times\mm{X}\times\mm{X}\times[0,1].$

  Тогда
    для всякой неограниченно возрастающей
    последовательности
  моментов времени $(\tau_n)_{n\in\mm{N}}\in \mm{T}^\mm{N}$
    найдется такое решение
    $(x^0,u^0,\lambda^0,\psi^0)\in\fr{Z}$
    всех соотношений принципа максимума
    \rref{sys_x}--\rref{dob},
     что   выполнено
    условие \rref{partlim_} для
    $A=A^{u^0}_{\xi^0},\ \xi^0=(0_\mm{X},\psi^0(0),\lambda^0).$
\end{zam}


\section{Формула Коши для сопряженной
переменной} \setcounter{equation}{0}

 В работе  \cite{kr_as} С.М.Асеевым,  А.В.Кряжимским
 был предложен и доказан
 вариант условия, в котором фактически в качестве матрицы $A$
 берется решение линеаризованной (вдоль оптимального решения)
  системы для уравнения \rref{sys_psi}.
 При этом для сопряженной переменной удается выписать в явном виде
 выражение, подобное формуле Коши  решения линейных уравнений.
 Это выражение обобщает
  (см.  \cite[\S 16]{kr_as}) целый ряд  условий трансверсальности,
  в частности, является более общей чем
 найденные    для линейных систем условия в \cite{aucl}.

 В отличие от большинства других условий трансверсальности, подход
 \cite{kr_as} хорош тем, что выделяется с каждым оптимальным
 управлением ровно одно решение принципа максимума. Таким образом,
 включение данного условия в соотношения принципа максимума  позволяет
 получить "полную"\ систему соотношений.

  Упростим предложение~\ref{3_} с таким~$A$
  для ослабления предположений
  \cite[теорема 2]{kr_asD},\cite[теорема 12.1]{kr_as}
  и их следствий.

\medskip

  Пусть некоторая пара
    $(x^0,u^0)\in C_{loc}(\mm{T},\mm{X})\times\fr{U}$
    оптимальна для задачи \rref{opt},
  определим вместе с ней решение задачи Коши:
 \begin{equation}
   \label{A0_def}
 \frac{d{A}(t)}{dt} =\frac{\partial f (t,x^0(t),u^0(t))}{\partial x} A(t),\
 A(0)=1_\mm{L}.
\end{equation}

  Подобным образом, для всякого  $\xi\in\mm{X}$ через $x_{\xi}$ обозначим решение
  \rref{sys_x} при начальном условии $x_{\xi}(0)=\xi\in\mm{X}$,
  введем также $A_{\xi}$ --- решение матричной задачи Коши
 $$\frac{d{A}_{\xi}(t)}{dt}=\frac{\partial
f(t,x_{\xi}(t),u^0(t))}{\partial x} A_{\xi}(t),\ A_{\xi}(0)=1_\mm{L}\
\ \forall \xi\in\mm{X},$$
 введем для всякого
  $T\in\mm{T}$ вектор
  \begin{equation}
   \label{omega2___}
   I_{\xi}(T)\rav\int_{[0,T]}
   \frac{\partial g(t,x_{\xi}(t),u^0(t))}{\partial x}\, A_\xi(t)
  \,dt.
\end{equation}

\begin{prp}
\label{aff}
    В условиях $\bf{(f),(g),(u),(\partial)}$
    пусть пара
    $(x^0,u^0)$
    оптимальна для задачи \rref{opt}.
 Пусть отображение $I_0$ ограничено и $\displaystyle \lim_{\xi\to 0} ||I_{\xi}-I_0||_C=0.$
 Пусть $I_*\in\mm{X}$ --- некоторый частичный предел
   $I_0(\tau)$ при $\tau\to \infty$.

 Тогда
    для
     $\lambda^0\rav 1/(1+||I_*||_\mm{X})>0$ и  $\psi\in C_{loc}(\mm{T},\mm{X}),$
     введенного по правилу: для всех\
      $T\in\mm{T}$
  \begin{equation}
   \label{klass_}
 \psi^0(T)\rav\frac{(I_*-I_0(T))A^{-1}(T)}{1+||I_*||_\mm{X}}=
 \lambda^0 \Big(I_*-\int_{[0,T]}
 \frac{\partial g(t,x^0(t),u^0(t))}{\partial x}\,A(t)\,dt\Big) A^{-1}(T),
\end{equation}
     четверка
    $(x^0,u^0,\lambda^0,\psi^0)$ удовлетворяет всем
    соотношениям принципа максимума
    \rref{sys_x}--\rref{dob}
    и условию  \rref{partlim_}.
\end{prp}
\doc
 Для всякого $\eta\in\td{\fr{U}}$ введём
 матричнозначную функцию $A^\eta$ --- решение матричного уравнения
\begin{equation}
   \label{Aeta}
 \dot{A}^\eta(t)=\int_{U(t)}\frac{\partial
f(t,\td{\ph}[\eta](t),u)}{\partial x}\,A^\eta(t)\, \eta(t)(du),\
A^\eta(0)=1_\mm{L},
\end{equation}
 теперь для всякого решения
$(x^\eta,\eta,\lambda^\eta,\psi^\eta)\in\td{\fr{Y}}$
 из
\rref{sys_psi_} следует
\begin{equation}
   \label{AetaPsi_1}
\frac{d}{dt}\big(\psi^\eta A^\eta\big)(t)=-\lambda^\eta
\int_{U(t)}\frac{\partial g(t,\td{\ph}[\eta](t),u)}{\partial x}
A^\eta(t)\eta(t)(du).
\end{equation}

  Превратим  систему
  \rref{AetaPsi_1},\rref{Aeta},\rref{sys_x},\rref{sys_lambda}
  в систему вида \rref{aa} подстановкой:
  $$E_p\rav\mm{X},\ E_q\rav\mm{L}\times\mm{X}\times\mm{R},\
   K_0\rav\{\psi\in\mm{X}\,|\, ||\psi||_\mm{X}\leq
    1\}\times\{1_\mm{L}\}\times\{0\}\times[0,1],\
    G_0\rav\mm{X}\times\mm{L}\times\mm{X}\times[0,1],$$
 $$   b(t,p,(q_1,q_2,q_3),u)\rav-q_3 \frac{\partial
 g(t,q_2,u)}{\partial x}  q_1,\ \  c(t,p,(q_1,q_2,q_3),u)\rav
 \Big( \frac{\partial
 f(t,q_2,u)}{\partial x}q_1,f(t,q_2,u),0\Big),
 $$

   Проверим частичную устойчивость по Ляпунову компоненты $p$ полученной системы
   $$
     \dot{p}=-r \frac{\partial
 g(t,z,u^0(t))}{\partial x} B , \dot{B}=\frac{\partial
 f(t,z,u^0(t))}{\partial x} B ,\dot{z}=f(t,z,u^0(t)),
 \dot{r}=0
     $$
   при
   $B(0)=1_\mm{L},z(0)=0,r(0)\in[0,1], p(0)\in\mm{D}.$
   Решая систему имеем:
 \begin{equation}
   \label{1943}
   r(t)=r(0),\ z(t)=x_{z(0)}(t),\ B(t)=A_{z(0)}(t)B(0),\
   p(t)=p(0)-r(0)I_{z(0)}(t)B(0).
\end{equation}
  Устойчивость по Ляпунову компоненты $p$
   на всем $K_0$ свелась к имеющим место по условию предложения непрерывности $I_{\xi}$ в точке $\xi=0$ и
   ограниченности $I_0.$
   Тогда по
    предложению  \ref{zvez} при некотором весе $w^0\in (Null)(u^0)$
    компонента  $p=A^\eta \psi^\eta$ устойчива при малых по $w^0$ возмущениях
   управления $u^0.$ В силу выбора $K_0$ это обеспечивает условие $\bf{(s'A)}$
   для так определенных $A^\eta$.

   По условию  $I_*$ --- частичный предел, потому найдется неограниченно возрастающая последовательность
 $(\tau_n)_{n\in\mm{N}}\in\mm{T}^\mm{N}$, для которой $I_0(\tau_n)\to
 I_*$, зафиксируем эту последовательность.
  Теперь по следствию \ref{s3_}
  для некоторой четверки $(x^0,u^0,\lambda^0,\psi^0)\in\fr{Z}$ выполнено
  \rref{partlim_}.

  Подставляя в \rref{1943} $z(0)=0, r(0)=\lambda^0,$   с учетом $A(0)=1_\mm{L}$,
  $p(0)=\psi^0(0)$, имеем
\begin{equation}
   \label{AetaPsi_2}
   p(T)=\big(\psi^0A\big)(T)=\psi^0(0)-\lambda^0I_0(T)\ \ \forall T\in\mm{T}.
\end{equation}
 Подставляя теперь уже $T=\tau_n$  и переходя  к пределу,
 из \rref{partlim_} следует $0=\psi^0(0)-\lambda^0I_*,$
 теперь из
 \rref{AetaPsi_2} и \rref{dob}
 соответственно, имеем
 $$\displaystyle \psi^0(T)A(T)=\lambda^0\big(I_*-I_0(T)\big),\ \lambda^0=\frac{1}{1+||I_*||_\mm{X}}> 0.$$
 Применяя матрицу, обратную к $A(T)$
 (решение задачи \rref{A0_def} невырождено), получаем \rref{klass_}.
 \bo

\medskip

 Отметим, что в том случае, когда $I_*$ не зависит от выбора
 последовательности $(\tau_n)_{n\in\mm{N}},$  автоматически
 выполнено более сильное условие трансверсальности
   \begin{equation}
   \label{lim}
   \lim_{t\to\infty} \psi(t)A(t)=0.
\end{equation}
  Более того, поскольку при фиксированных $(x^0,u^0,\lambda)$ решения \rref{AetaPsi_1}  отличаются друг от
  друга на константу, то для всех остальных четверок
  $(x^0,u^0,\lambda,\psi)\in\fr{Y}$ произведения $\psi A$
  также будут стремиться при $t\to\infty$ к конечному пределу, но уже
  не равному нулю. Отсюда каждой тройке $(x^0,u^0,\lambda)$
  соответствует не более одного $\psi^0$, для которого выполнено
  \rref{sys_x}--\rref{maxH},
  \rref{lim};
  теперь из \rref{dob} однозначно восстанавливается и $\lambda^0$.
  Итак, существует единственное решение
  $(x^0,u^0,\lambda,\psi)\in\fr{Z}$,
  удовлетворяющее условию \rref{lim}. Таким образом показана:

\begin{teo}
\label{aff_}
    В условиях $\bf{(f),(g),(u),(\partial)}$
    пусть пара
    $(x^0,u^0)$
    оптимальна для задачи \rref{opt},
    кроме того определен, конечен и непрерывно зависит от $\xi$ в некоторой окрестности $0_\mm{X}$
    предел
   \begin{equation}
   \label{omega2}
   \Lambda_{\xi}\rav\lim_{T\to\infty}
  \int_{[0,T\ra}
  \frac{\partial g(t,x_{\xi}(t),u^0(t))}{\partial x}\,A_\xi(t) dt
\end{equation}

 Тогда существует единственная четверка $(x^0,u^0,\lambda^0,\psi^0)\in\fr{Z}$,
  удовлетворяющая всем соотношениям
 принципа максимума     \rref{sys_x}--\rref{dob}
и
    условию трансверсальности \rref{lim}. При этом
  \begin{equation}
   \label{klass}
 \lambda^0\rav
 \frac{1}{1+||\Lambda_0||_\mm{X}}>0,\ \ \
 \psi^0(T)\rav
 \int_{[T,\infty\ra}
 \frac{\partial g(t,x^0(t),u^0(t))}{\partial x}  A(t) dt \frac{A^{-1}(T)}{1+||\Lambda_0||_\mm{X}}\ \ \forall T\in\mm{T}.
\end{equation}
\end{teo}

В
 \cite{kr_asD},\cite{kr_as} рассматривалась задача с функцией
 $f$,  не зависящей от времени и функцией $g$ вида $e^{-rt}G(x,u)$;  накладываемые в
  \cite[теорема~2]{kr_asD},\cite[теорема
 12.1]{kr_as} условия (доминирования дисконтирующего множителя) в терминах условия
  \rref{omega2} могут быть сведены к виду:
   при некоторых $\alpha,\beta\in\mm{R}_{>0}$
  для всех допустимых управлений $u$, соответствующих им траекторий $x$ и матричных
  экспонент $A$
  выполнено:
 $$   ||A(t)||\,
  \Big|\Big|\frac{\partial g(t,x(t),u(t))}{\partial x}\Big|\Big| \leq \beta e^{-\alpha t}\
  \ \forall t\in\mm{T}.$$
 Это условие, в свою очередь, может быть проверено подсчетом
 показателей Ляпунова системы принципа максимума, смотрите \cite[\S 12]{kr_as}.

 Отметим также, что формула \rref{klass} показана также для биафинной
 управляемой
 системы при монотонной $\frac{\partial g}{\partial x}$  (\cite[теорема~1]{kr_asD},\cite[теорема
 11.1]{kr_as}). Этот результат,
по-видимому
  напрямую из показанной в данной работе теоремы \ref{aff_} не следует.

 В работе \cite[Example 2]{ss} рассматривается следующий
\begin{example}
$$\dot{x}=-u,\ ,\ x(0)\in\mm{R}_{>0},\  u\in[c,d]\subset\mm{R}_{>0},\ \ \int_{\mm{T}}e^{-r(t)}(s(x)+u)dt\to \max, $$
\end{example}
  где  
 скалярные функции $s,r$ непрерывно дифференцируемы,
 $|s'(x)|$ растет не
 быстрее некоторого полинома и существует $\lim_{t\to\infty}\dot{r}(t)>0$.
  Как утверждается в \cite{ss}
 при $r(t)=t,\ s(x)=x^3$ стандартные условия
 на бесконечности отсутствуют, зато удаётся показать
\begin{equation}
   \label{2}
   \lim_{t\to+\infty}\ct{H}(x^0(t),t,u^0(t),\lambda^0,\psi^0(t))=0.
\end{equation}
 Но здесь  $A(t)\equiv 1_{\mm{L}},$ интеграл
 $\int_{\mm{T}}e^{-r(t)}s'(x(t))dt$ равномерно по $x$ сходится и выполнено  \rref{omega2}.
 Тогда из теоремы \ref{aff_} имеем $\lambda^0>0,$
 $\psi^0(t)=\lambda^0\int_{t}^\infty e^{-r(t)}s'(x^0(t))dt,$
  а  следовательно \rref{partlim} и  \rref{2}.

 \noindent Хлопин Дмитрий Валерьевич\hfill Поступила  \\
Ин-т математики и механики УрО РАН\\
e-mail: khlopin@imm.uran.ru\\
канд. физ.-мат.\ наук\\
ст.\ науч. сотрудник\\[1ex]
\end{document}